\documentclass[a4]{article}
\usepackage{graphicx}
\newtheorem{thm}{Theorem}

\def\be{\begin{equation}}
\def\ee{\end{equation}}
\newcommand{\kk}[2]{\frac{#1}{#2}}

\newcommand{\ff}[1]{{\bf #1}}

\def\Q{\Omega}

\def\lam{\lambda}
\def\ra{\rightarrow}
\def\x{\ff{x}}

\newcommand{\brak}[1]{\left\{  \begin{array}{lllll} #1 \end{array} \right. }
\begin{document}

\title{Global Convergence Analysis of the Flower Pollination Algorithm:
A Discrete-Time Markov Chain Approach }

\author{Xingshi He$^1$, Xin-She Yang$^2$, Mehmet Karamanoglu$^2$, Yuxin Zhao$^3$ \\[10pt] 
1) College of Science, Xi'an Polytechnic University, Xi'an, P. R. China  \\
2) School of Science and Technology, Middlesex University, London NW4 4BT, UK\\
3) College of Automation, Harbin Engineering University, Harbin, P. R. China  \\
}

\date{}

\maketitle

\begin{abstract}
Flower pollination algorithm is a recent metaheuristic algorithm for solving
nonlinear global optimization problems. The algorithm has also been
extended to solve multiobjective optimization
with promising results. In this work, we analyze this algorithm mathematically
and prove its convergence properties by using Markov chain theory.
By constructing the appropriate transition probability for a population
of flower pollen and using the homogeneity property, it can be shown
that the constructed stochastic sequences can converge to the optimal set.
Under the two proper conditions for convergence, it is proved
that the simplified flower pollination algorithm can indeed satisfy these
convergence conditions and thus the global
convergence of this algorithm can be guaranteed.
Numerical experiments are used to demonstrate that the  flower pollination
algorithm can converge quickly in practice and can thus achieve global optimality efficiently.  \\[10pt]

\noindent {\bf Citation Detail:} X.S. He, X.-S. Yang, M. Karamanoglu, Y.X. Zhao, Global convergence analysis of the flower pollination algorithm: A discrete-time Markov chain approach, {\it Procedia Computer Science}, vol. 108, 1354-1364 (2017).
https://doi.org/10.1016/j.procs.2017.05.020

\end{abstract}

\section{Introduction}

Computational intelligence and optimization have become increasingly important in
many applications, partly due to the explosion of data volumes driven by the Internet
and social media, and partly due to the more stringent design requirements.
In recent years, bio-inspired optimization algorithms have gained some popularity \cite{Kennedy,Yang2014}.
Many new optimization algorithms are based on the so-called swarm intelligence with diverse
characteristics in mimicking natural systems. Consequently, different algorithms may
have different features and thus may behave differently, even with different efficiencies.
However, it still lacks in-depth understanding why these algorithms work well and exactly under what conditions.

In fact, there is a significant gap between theory and practice. Most metaheuristic algorithms have
successful applications in practice, but their mathematical
analysis lacks far behind. In fact, apart from a few limited results about the convergence
and stability concerning particle swarm optimization, genetic algorithms, simulated annealing and others
\cite{Clerc,Jiang,Ren,YangRW}, many algorithms do not have theoretical analysis.
Therefore, we may know they can work well in practice, but we rarely understand why they work
and how to improve them with a good understanding of their working mechanisms.

Among most recent, bio-inspired algorithms, flower pollination algorithm (FPA), or flower algorithm (FA) for simplicity,
has demonstrated very good efficiency in solving both single objective optimization and multi-objective optimization problems \cite{Yang2012FPA,Yangetal}.
This algorithm mimics the main characteristics of the pollination
process of flowering plants, which leads to
both local and global search capabilities. As this algorithm is very new, there is no mathematical analysis yet.

The main purpose of this paper is to analyze the flower algorithm mathematically and try to
prove its convergence properties. Therefore, this paper is organized as follows.
In Section 2, the flower algorithm will be outlined briefly, followed by some
simplifications so as to be used for the detailed mathematical analysis in Section 3 and Section 4.
Then, in Section 5, some numerical benchmarks will be used to demonstrate the
main characteristics of the convergence behaviour of the flower algorithm. Finally, conclusions
will be drawn briefly in Section 6.

\section{Flower Pollination Algorithm and Applications}

\subsection{Flower Algorithm}

Flower pollination algorithm (FPA), or flower algorithm,
was developed by Xin-She Yang in 2012 \cite{Yang2012FPA},
inspired by the flow pollination process of flowering plants.
The flower pollination algorithm has then been extended to deal with multiobjective optimization \cite{Yang2013,Yangetal}. The diversity of flowering plants are amazing, and it  is estimated that there are over a quarter of a million types of flowering plants in Nature and that
about 80\% of all plant species are flowering species.  Flower pollination is
typically associated with the transfer of pollen, and such transfer
is often linked with pollinators such as insects, birds, bats and other animals.
Pollination can take two major forms: abiotic and biotic. About 90\%
of flowering plants belong to biotic pollination. That is, pollen is
transferred by a pollinator such as insects, bats and animals.
In fact, some flowers and insects have co-evolved into a very specialized
flower-pollinator partnership called flower constancy \cite{Waser}.
For example, hummingbirds are a good example
for flower constancy in pollination. Such flower constancy may have evolutionary advantages
because this will maximize the transfer of flower pollen to the same or conspecific plants,
and thus maximizing the reproduction of the same flower species. For the pollinators in the
flower constancy partnership,
they can minimize their efforts for searching for new flower patches and thus
with a higher probability of nectar rewards from the same flower species.

Pollination can be achieved by self-pollination or cross-pollination.
Self-pollination tends to be local and often occurs when there is no reliable pollinator available.
On the other hand, biotic cross-pollination may occur at long distance, and the
pollinators such as bees, bats, birds and flies can fly a long distance,
thus they can be considered as the global pollination.

For simplicity in describing the flower algorithm, the following four rules can be summarized as follows \cite{Yang2012FPA,Yangetal}:

\begin{enumerate}

\item Biotic and cross-pollination can be considered as a process of
global pollination, and pollen-carrying pollinators move in a way which
obeys L\'evy flights (Rule 1).

\item For local pollination, abiotic and self-pollination can be used (Rule 2).

\item Pollinators such as insects can develop flower constancy, which is equivalent to
a reproduction probability that is proportional to the similarity of two flowers involved (Rule 3).

\item The interaction or switching of local pollination and global pollination can
be controlled by a switch probability $p \in [0,1]$, with a slight bias towards local
pollination (Rule 4).

\end{enumerate}
In order to formulate the updating formulae in the FPA,  we have to convert the above rules into
updating equations. For example, in the global pollination step,
flower pollen gametes are carried by pollinators such as insects,
and pollen can travel over a long distance because insects can often fly and travel in a much
longer range. Therefore, Rule 1 and flower constancy can be  represented mathematically as
\be \x_i^{t+1}=\x_i^t + \gamma L(\lam) (\ff{g}_*-\x_i^t ), \ee
where $\x_i^t$ is the pollen $i$ or solution vector $\x_i$ at iteration $t$, and $\ff{g}_*$ is
the current best solution found among all solutions at the current generation/iteration.
Here $\gamma$ is the parameter that corresponds to the strength of the pollination,
which essentially is also a step size.
Since insects may move over a long distance with various distance steps, we can use a L\'evy flight
to mimic this characteristic efficiently. That is, we draw $L$ from a L\'evy distribution \cite{Yangetal,Pav}
 \be L \sim \frac{\lam \Gamma(\lam) \sin (\pi \lam/2)}{\pi} \frac{1}{s^{1+\lam}}, \quad (s \gg s_0>0). \ee
Here $\Gamma(\lam$) is the standard gamma function, and this distribution is valid for
large steps $s>0$. Though in theory the critical size $s_0$ should be sufficiently large, $s_0=0.1$ or even $s_0=0.01$ can be used in practice. Here, the notation `$\sim$' means to draw random numbers that obey the distribution on the right-hand side.

For the local pollination, both Rule 2 and Rule 3 can be represented as
\be \x_i^{t+1}=\x_i^t + U (\x_j^t -\x_k^t), \ee
where $\x_j^t$ and $\x_k^t$ are pollen from different flowers of the same plant species.
This essentially mimics the flower constancy in a limited neighborhood. Mathematically,
if $\x_j^t$ and $\x_k^t$ comes from the same species
or selected from the same population, this equivalently becomes a local random walk if we draw
$U$ from a uniform distribution in [0,1].

In principle, flower pollination activities can occur at all scales, both local and global.
But in reality,  adjacent flower patches or flowers in the not-so-far-away neighborhood
are more likely to be pollinated by local flower pollen than those far away.
In order to mimic this feature, we can effectively use a
switch probability (Rule 4) or proximity probability $p$ to switch
between common global pollination to intensive local pollination.
To start with, we can use a naive value of $p=0.5$ as an initially value.
A parametric study showed that $p=0.8$ may work better for most applications.
Preliminary studies suggest that the flower algorithm is very efficient, and has been
extended to multi-objective optimization \cite{Yang2013,Yangetal}.

It is worth pointing out that parameter tuning may be needed in all algorithms, and ideally a self-tuning framework can be used \cite{Yang2013ST}. However, in our analysis of convergence, we assume that the parameter values are fixed, though such parameter values can be within a range.
In addition, the representations of the solution vectors in the algorithm are simply vectors,
not in any complicated forms such as quaternion representations \cite{FisterQuat}.

\subsection{Applications}

Since the development of the basic flower pollination algorithm (FPA), there are a wide range of diverse applications of this algorithm with more than 500 research papers published so far in the literature. For example, a brief review by Chiroma et al. identified some of the earlier applications \cite{Chiroma}. Therefore, it is not possible to review even a small fraction of the latest developments. Here, we only highlight a few recent papers. For example, Dubey et al. presented a hybrid FPA variant for
solving multi-objective economic dispatch problems \cite{Dubey,Dubey2}, while Alam et al. carried out  photovoltaic parameter estimation using FPA \cite{Alam}. Structure optimization has also been investigated using FPA \cite{Bekdas}, and feature selection has been done using a clonal FPA by
Sayed et al. \cite{Sayed}. A modified FPA for global optimization has been proposed by Nabil \cite{Nabil}.

In addition, Velamuri et al. used FPA to optimize economic load dispatch \cite{Velam}, while Rodrigues et al. developed a binary flower pollination algorithm to do EEG identification.
Furthermore, Zhou et al. introduced an elite opposition-based FPA \cite{Zhou} and Mahdad et al. presented an adaptive FPA to solve optimal power flow problems \cite{Mahdad}, while
Abdelaziz et al. solved placement problems in distribution systems using FPA \cite{Abdel}. New variants of FPA are still emerging \cite{Salgotra}.

Obviously, there are other important applications, but here we will focus on the mathematical analysis of the basic FPA. Therefore, we will start with the simplified version of FPA.

\subsection{Simplified Flower Algorithm}
As there are two branches in the updating formulas, the local search step only contributes mainly
to local refinements, while the main mobility or exploration is carried out by
the global search step. In order to simplify the analysis and also to emphasize the global search
capability, we now use a simplified version of the flower algorithm. That is, we use only the global
branch with a random number $r \in [0,1]$, compared with a discovery/switching probability $p$.
Now we have
\be \brak{\x_i^{(t+1)} \leftarrow x_i^{(t)} & \textrm{ if } r>p, \\ \\
\x_i^{(t+1)} \leftarrow  \x_i^{(t)} + \gamma G & \textrm{ if } r<p, } \ee
where $G(\gamma, \ff{g}_*, \x_i^t)=L(\lam) (\ff{g}_*-\x_i^t)$.

As the flower pollination algorithm is a stochastic search algorithm, we can summarize the simplified version
as the following key steps:

\begin{itemize}

\item[Step 1.] Randomly generate an initial population of $n$ pollen agents
at the positions, $\ff{X}=\{\x_1^0, \x_2^0, ..., \x_n^0 \}$,
then evaluate their objective values so as to find the initial best $\ff{g}_t^0$.

\item[Step 2.] Update the new solutions/positions by
\be \x_i^{(t+1)} =\x_i^{(t)} + \gamma G.  \label{equ-100} \ee

\item[Step 3.] Draw a random number $r$ from a uniform distribution $[0,1]$. Update $\x_i^{(t+1)}$ if
$r>p$. Then, evaluate the new solutions so as to find the new, global best $\ff{g}_t^*$ at pseudo time/iteration $t$.

\item[Step 4.] If the stopping criterion is met, then $\ff{g}_t^*$ is the best global solution found so far.
Otherwise, return to Step 2 and continue.

\end{itemize}

\section{Convergence Analysis}

\subsection{Gap Between Theory and Practice}

There is a significant gap between theory and practice in bio-inspired computing.
Nature-inspired metaheuristic algorithms work almost magically in practice, but it is not well
understood why these algorithms work. For example, except for a few cases such as genetic
algorithms, simulated annealing and particle swarm optimization, there are not many good
results concerning the convergence analysis and stability of metaheuristic algorithms. The lack
of theoretical understanding may lead to slow progress or even resistance to the wider
applications of metaheuristics.

There are three main methods for theoretical analysis of algorithms, and they are: complexity
theory, dynamical systems and Markov chains. On the one hand, metaheuristic algorithms tend
to have low algorithm complexity, but they can solve highly complex problems.
On the other hand, the convergence analysis typically use dynamic systems
and statistical methods as well as Markov chains. For example, particle swarm optimization was
analysed by Clerc and Kennedy \cite{Clerc} using simple dynamic systems, while genetic algorithms
was analysed intensively in a few theoretical studies \cite{Aytug,Green,Gut,Villa}.

For a genetic algorithm with a given mutation rate ($\mu$), string length ($L$) and population size ($n$), the number of iterations in genetic algorithm can be estimated by
\be t \le \Big\lceil \frac{\ln (1-p)}{\ln \Big(1-\min[(1-\mu)^{n L}, \mu^{n L}] \}
\Big)} \Big\rceil, \ee
where $\lceil u \rceil$ means taking the maximum integer value of $u$, and $p$ is a function of
$\mu$, $L$ and $n$ \cite{Gut,Villa}. However, for other bio-inspired algorithms, especially new algorithms,
theoretical understanding lacks behind, and thus there is a strong need for further studies
in this area. There is no doubt that any new understanding will provide greater insight into the
working mechanism of metaheuristic algorithms.

\subsection{Convergence Criteria in Stochastic Search}

For an optimization problem $<\!\Q,f\!>$, a stochastic search algorithm $A$, the $k$th iteration will produce a new solution
\be \x_{k+1}=A(\x_k, \xi), \ee
where $\Q$ is the feasible solution space, and $f$ is the objective function.
Here, $\xi$ is the visited solutions of algorithm $A$ during the iterative process.

In the Lebesgue measure space, the infimum of the search can be defined as
\be \phi=\inf \Big(t: \nu[x \in \Q\Big| f(x)<t]>0 \Big), \label{equ-200} \ee
where $\nu[X]$ denotes the Lebesque measure on the set $X$. Here Eq.(\ref{equ-200}) represents
the non-empty set in the search space, and the region for optimal solutions can be defined as
\be R_{\epsilon, M}=\brak{\{ x \in \Q | f(x) < \phi +\epsilon \}  & \textrm{ if } \phi \textrm{ is finite}, \\ \\
\{ x \in \Q | f(x) < -C \} & \textrm{ if } \phi=-\infty, } \ee
where $\epsilon>0$ and $C$ is a sufficiently large positive number.  If any point in $R_{\epsilon, M}$ is
found during the iteration, we can say the algorithm $A$ has found the globally optimal solution
or its best approximation.

In order to analyze the convergence of an algorithm, let us first state the
conditions for convergence \cite{Jiang,Wang}:

\begin{itemize}

\item Condition 1. If $f(A(x,\xi)) \le f(x)$ and $\xi \in \Q$, then $f(A(x,\xi)) \le f(\xi)$.

\item Condition 2. For $\forall B \in \Q$ subject to $\nu(B)>0$, $$\prod_{k=0}^{\infty} (1-u_k(B))=0, $$
where $u_k(B)$ is the probability measure on $B$ of $k$th iteration of the algorithm $A$.

It is worth pointing out that we focus on the minimization problems in our discussions.
\end{itemize}

{\bf Lemma 1.}
{\it The global convergence of an algorithm.} If $f$ is measurable and the feasible solution
space $\Q$ is a measurable subset on $\Re^n$, algorithm $A$ satisfies the above two conditions with
the search sequence $\{x_k\}_{k=0}^{\infty}$, then
\be \lim_{k \rightarrow \infty} P (x_k \in R_{\epsilon, M})=1. \ee
That is, algorithm $A$ will converge globally \cite{Jiang,Wang}. Here $P(x_k \in R_{\epsilon,M})$ is the probability
measure of the $k$th solution on $R_{\epsilon, M}$  at the $k$th iteration.

\section{Markov Chains and Convergence Analysis}

\subsection{Definitions}

{\bf Definition 1.} {\it The state and state space}. The positions of pollen and its global best solution $g$ in the search history forms the states of flower pollen: $y=(x,g)$, where $x, g \in \Q$ and $f(g) \le f(x)$ (minimization problems).
The set of all the possible states form the state space, denoted by
\be Y=\{y=(x,g)|x, g \in \Q, f(g) \le f(x) \}. \ee

{\bf Definition 2.} {\it The states and state space
of the pollen group/population}. The states of all $n$ solutions
form the states of the group, denoted by $q=(y_1, y_2, ..., y_n)$. All the states of all the pollen form a state space for the group, denoted by
\be Q=\{q=(y_1, y_2, ..., y_n), y_i \in Y, 1 \le i \le n \}. \ee
Obviously, $Q$ contains the historical global best solution $g^*$ for the whole population
 and all individual best solutions $g_i (1 \le i \le n)$ in history.
 In addition, the global best solution of the whole population is the best
 among all $g_i$, so that $f(g^*)=\min(f(g_i)), \; 1 \le i \le n$.

 \subsection{Markov Chain Model for Flower Algorithm}

{\bf Definition 3.} {\it The state transition for pollen positions}. For $\forall y_1=(x_1, g_1) \in Y,
\forall y_2=(x_2, g_2) \in Y$, the state transition from $y_1$ to $y_2$ can be denoted by
\be T_y(y_1)=y_2. \ee

\begin{thm}
The transition probability from state $y_1$ to $y_2$ in the flower algorithm is
\[ P(T_y(y_1)=y_2)=P(x_1 \ra x_1') P(g_1 \ra g_1') \]
\be \times P(x_1' \ra x_2) P(g_1' \ra g_2), \ee
\end{thm}
where $P(x_1 \ra x_1')$ is the transition probability at Step 2 in the flower algorithm,
and $P(g_1 \ra g_1')$ is the transition probability for the historical global best at this step.
$P(x_1' \ra x_2)$ is the transition probability at Step 3, while $P(g_1' \ra g_2)$ is
the transition probability of the historical global best.

{\it Proof:} In the simplified flower algorithm,  the state transition from $y_1$ to $y_2$ only
has one middle transition state $(x_1', g_1')$, which means
that $x_1 \ra x_1', g_1 \ra g_1', x_1' \ra x_2$ and $g_1' \ra g_2$ are valid simultaneously. Then, the probability for $P(T_y(y_1)=y_2)$ is
\[ P(T_y(y_1)=y_2)=P(x_1 \ra x_1') P(g_1 \ra g_1') \]
\be \times  P(x_1' \ra x_2) P(g_1' \ra g_2). \label{equ-300} \ee

From Eq. (\ref{equ-100}), the transition probability for $x_1 \ra x_1'$ is
\be P(x_1 \ra x_1')=\brak{\frac{1}{|g-x_1|} & \textrm{ if } x_1' \in [x_1, x_1+(x_1-g)], \\ \\
0 & \textrm{ if } x_1' \notin [x_1, x_1+(x_1-g)]. } \label{equ-400} \ee
Since $x$ and $g$ are higher-dimensional vectors, the mathematical operations here
should be interpreted as vector operations, while the $|\cdot|$ means the volume of the hypercube.

The transition probability of the historical best solution is
\be P(g_1 \ra g_1') =\brak{1 & f(x_1') \le f(g_1), \\  \\ 0 & f(x_1') > f(g_1). } \label{equ-400} \ee
From Step 3 in the simplified
flower algorithm, we know that a random number $r \in [0,1]$ is compared
with the discovery probability $p=4/5$. If $r<p$, then the position/solution of pollen can be changed randomly;
otherwise, it remains unchanged. Therefore,
the transition probability for $x_1' \ra x_2$ is
\be P(x_1' \ra x_2) =\brak{ 1-p, &  \textrm{ if } \; r>p, \\ \\ p,   & \textrm{ if } \; r \le p, }
=\brak{\frac{4}{5} & \textrm{ if }\; r> p, \\ \\ \frac{1}{5} & \textrm{ if } \; r \le p. } \ee
The transition probability for the historical best solution is
\be P(g_1' \ra g_2) =\brak{1 & \textrm{ if } f(x_2) \le f(g_1), \\ 0 & \textrm{ if } f(x_2) > f(g_1). } \label{equ-700} \ee
{\bf Definition 4.} {\it The group transition probability in flower algorithm}. The group transition
probability can be defined as $T_q(q_i)=q_j$ for
$\forall q_i=(y_{i1}, y_{i2}$, $..., y_{in}) \in \Q$ and $\forall q_j=(y_{j1}, y_{j2}, ..., y_{jn}) \in \Q$.
\begin{thm}
In the simplified flower algorithm, the group transition probability from $q_i$ to $q_j$ in one step is
\be P(T_q(q_i)=q_j) =\prod_{k=1}^n P(T_y(y_{ik})=y_{jk}).  \label{equ-800} \ee
\end{thm}
{\it Proof:} If the group states can be transferred from $q_i$ to $q_j$ in one step, then all the states
will be transferred simultaneously. That is, $T_y(y_{i1}=y_{j1}, T_y(y_{i2})=y_{j2}$, ..., $T_y(y_{in})=y_{jn}$, and
the group transition probability can be written as the joint probability
\[ P(T_q(q_i)=q_j)=P(T_y(y_{i1})=y_{j1}) P (T_y(y_{i2})=y_{j2})
\] \be \cdots P(T_y(y_{in})=y_{jn})  =\prod_{k=1}^n P(T_y(y_{ik})=y_{jk}). \ee

\begin{thm}
The state sequence $\{q(t); t \ge 0\}$ in the flower algorithm is a finite homogeneous Markov chain.
\end{thm}

{\it Proof:} First, let us assume that all search spaces for a stochastic algorithm are finite.
Then, $x$ and $g$ in any pollen state $y=(x,g)$ are also finite, so that the state space for
flower pollen are finite. Since the group state $q=(y_1, y_2, ..., y_n)$ consists of $n$ positions
where $n$ is positive and finite, so group states $q$ are also finite.

From the previous theorems, we know that the group transition probability
\be P(T_q(q(t-1))=q(t), \ee
for $\forall q(t-1) \in Q$ and $\forall q(t) \in Q$ is the group transition probability
$P(T_y(y_i(t-1))=y_i(t))$ for $1 \le i \le n$. From Eq. (\ref{equ-300}), we have the transition
probability for any pollen is
\[ P(T_y(y(t-1))=y(t))  =P(x(t-1) \ra x'(t-1)) \]
\[ \times P (g(t-1) \ra g'(t-1))   P(x'(t-1) \ra x(t)) \]
\be \times P (g'(t-1) \ra g(t)), \ee
where $P(x(t-1) \ra x'(t-1))$, $P(g(t-1) \ra g'(t-1))$, $P(x'(t-1) \ra x(t))$
and $P(g'(t-1) \ra g(t))$ are all only depend on $x$ and $g$ at $t-1$.
Therefore, $P(T_q(q(t-1))=q(t))$ also only depends on the states $y_i(t-1), 1 \le i \le n$ at time $t-1$.
Consequently, the group state sequence $\{q(t); t \ge 0 \}$ has the property of a Markov chain.

Finally, $P(x(t-1) \ra x'(t-1))$, $P(g(t-1) \ra g'(t-1))$, $P(x'(t-1) \ra x(t))$
and $P(g'(t-1) \ra g(t))$ are all independent of $t$, so is $P(T_y(y(t-1))=y(t))$.
Thus, $P(T_q(q(t-1))=q(t)$ is also independent of $t$, which implies that
this state sequence is also homogeneous. In summary,
the group state sequence $\{q(t); t\ge 0 \}$ is a finite, homogeneous Markov chain.

\subsection{Global Convergence of the Flower Algorithm}

{\bf Definition 5.} For the globally optimal solution $g_b$ for an optimization (or minimization)  problem
$<\!\Q,f\!>$, the optimal state set is defined as $R=\{y=(x,g)|f(g)=f(g_b), y \in Y\}$.

\begin{thm}
Given the position state sequence $\{y(t); t\ge 0\}$ in the flower algorithm, the state set $R$
of the optimal solutions corresponding to optimal solutions form a closed set on $Y$.
\end{thm}
{\bf Proof:} For $\forall y_i \in R, \forall y_j \notin R$, the probability for
$T_y(y_j)=y_i$ is $$P(T_y(y_j)=y_i)=P(x_j \ra x_i') P(g_j \ra g_j') P(x'_j \ra x_i) P (g'_j \ra g_j). $$
Since for $\forall y_i \in R$ and $\forall y_j \notin R$, it holds that
$f(g_i) \ge f(g_j)=f(g_b)=\inf(f(a)), a \in \Q$.

From Eqs. (\ref{equ-400}) and (\ref{equ-700}), we have
$P(g_j \ra g_j') P(g_j' \ra g_i)=0$, which leads to
$P(T_y(y_j)=y_i)=0$. This condition implies that $R$ is closed on $Y$.

{\bf Definition 6.} For the globally optimal solution $g_b$ to an optimization
problem $<\!\Q,f\!>$, the optimal group state set can be defined as
\be H=\{q=(y_1, y_2, ..., y_n)|\exists y_i \in R, 1 \le i \le n\}. \ee

\begin{thm}
Given the group state sequence $\{q(t); t \ge 0\}$ in the flower algorithm,
the optimal group state set $H$ is closed on the group state space $Q$.
\end{thm}

{\it Proof:} From Eq. (\ref{equ-800}), the probability
\be P(T_q(q_j)=q_i)=\prod_{k=1}^n P(T_y(y_{jk})=y_{ik}), \ee
for $\forall q_i \in H, \forall q_j \in H$ and $T_q(q_j)=q_i$.
Since $\forall q_i \in H$ and $\forall q_j \notin H$, in order to satisfy
$T_q(q_j)=q_i$, there exists at least one position that will transfer
from the inside of $R$ to the outside of $R$. That is,
$\exists T_y(y_{jk})=y_{ik}, y_{jk} \in R, y_{ik} \notin R,
1 \le k \le n$. From the previous theorem, we know that $R$ is closed on $Y$, which
means that $P(T_y(y_{jk})=y_{ik})=0$. Therefore,
$$ P(T_q(q_j)=q_i)=\prod_{k=1}^n P(T_y(y_{jk})=y_{ik})=0. $$
From the definition of a closed set, we can conclude that
the optimal set $H$ is also closed on $Q$.

\begin{thm}
In the group state space $Q$ for flower pollen, there does not exist
a non-empty closed set $B$ so that $B \cap H =\emptyset$.
\end{thm}

{\it Proof:} {\it Reductio ad absurdum.} Assuming that there exists a closed set $B$ so
that $B \cap H =\emptyset$ and that $f(g_j) > f(g_b)$ for
$q_i=(g_b, g_b, ..., g_b) \in H$ and $\forall q_j=(y_{j1}, y_{j2}, ..., y_{jn}) \in B$, then Eq. (\ref{equ-800}) implies that
\be P(T_q(q_j)=q_i)=\prod_{k=1}^n P(T_y(y_{jk})=y_{ik}). \ee
For each $P(T_y(y_j)=y_i)$, it holds that
$P(T_y(y_j)=y_i)=P(x_j \ra x_j') P(g_j \ra g_j') P(x_j' \ra x_i) P(g_j' \ra g_i)$.
Since $P(g_j' \ra g_i)=1, P(g_j \ra g_j'), P(x_j \ra x_j') P(x_j' \ra x_i) >0$, then
$P(T_y(y_j)=y_i) \ne 0$,  implying that $B$ is not closed, which contradicts
the assumption. Therefore, there exists no non-empty closed set outside $H$ in $Q$.

With the above definitions and results, it is straightforward to prove the following lemma:

{\bf Lemma 2.} Assuming that a Markov chain has a non-empty set $C$ and there does not exist
a non-empty closed set $D$ so that $C \cap D=\emptyset$, then
$\lim_{n \ra \infty} P(x_n =j) =\pi_j$ only if $j \in C$, and $\lim_{n \ra \infty} P(x_n =j) =0$ only
if $j \notin C$.

In addition, we have also have the following theorem:
\begin{thm}
When the number of iteration approaches infinity, the group state
sequence will converge to the optimal state/solution set $H$.
\end{thm}

{\it Proof:} Using the previous two theorems and Lemma 2, it is straightforward to prove this theorem.

Now it is ready to state the global convergence theorem.
\begin{thm}
The flower algorithm has guaranteed global convergence.
\end{thm}
{\it Proof:} First, the iteration process in the flower algorithm always keeps/updates the
current the global best solution for the whole population, which ensures that it satisfies the
the first convergence condition as outlined in the earlier section. From the previous theorem, the group state sequence will
converge towards the optimal set after a sufficiently large number of iterations or infinity.
Therefore, the probability of not finding the globally optimal solution is $0$, which
satisfies the second convergence condition. Consequently, the flower algorithm has guaranteed
global convergence towards its global optimality.

\section{Global Convergence and Numerical Experiments}

Many optimization algorithms are local search algorithms, though most metaheuristic algorithms
tend  to be suitable for global optimization. For multimodal objectives with many
local modes, many algorithms may be trapped in a local optimum.
As we have shown that the flower algorithm has good global convergence property, it can be
particularly suitable for global optimization. In order to show that  the flower algorithm indeed has
good convergence for various functions, we have chosen 5 different functions with
diverse modes and properties.

The Ackley function \cite{Ackley} can be written as
\be f(\x)=-20 \exp\Big[-\kk{1}{5} \sqrt{\kk{1}{d} \sum_{i=1}^d x_i^2}\Big]
 - \exp\Big[\kk{1}{d} \sum_{i=1}^d \cos (2 \pi x_i)\Big]
+20 +e, \ee
which has the global minimum $f_*=0$ at $(0,0,...,0)$.

The simplest of De Jong's functions is the so-called sphere function
\be f(\x) =\sum_{i=1}^n x_i^2, \quad
-5.12 \le x_i \le 5.12, \ee
whose global minimum is obviously $f_*=0$ at $(0,0,...,0)$. This function is unimodal and convex.

Rosenbrock's function \be f(\x) = \sum_{i=1}^{d-1} \Big[ (x_i-1)^2 + 100 (x_{i+1}-x_i^2)^2 \Big], \ee
whose global minimum $f_*=0$ occurs at $\x_*=(1,1,...,1)$ in the domain
$-5 \le x_i \le 5$ where $i=1,2,...,d$. In the 2D case, it is often written as
\be f(x,y)=(x-1)^2 + 100 (y-x^2)^2, \ee
which is often referred to as the banana function.

Yang's forest-like function
\be f(\x)=\Big( \sum_{i=1}^d |x_i| \Big) \exp\Big[- \sum_{i=1}^d \sin (x_i^2) \Big],
\quad -2 \pi \le x_i \le 2 \pi, \ee
has the global minimum $f_*=0$ at $(0,0,...,0)$, though the objective at this point is non-smooth.

Zakharov's function
\be f(\x)=\sum_{i=1}^d x_i^2 + \Big(\sum_{i=1}^d \frac{i x_i}{2} \Big)^2
+\Big(\sum_{i=1}^d \frac{i x_i}{2} \Big)^4, \ee
has its global minimum $f(\x_*)=0$ at $\x_*=(0,0, ..., 0)$ in the domain
$-5 \le x_i \le 5$.

\begin{figure}
\centerline{\includegraphics[height=2.5in,width=4in]{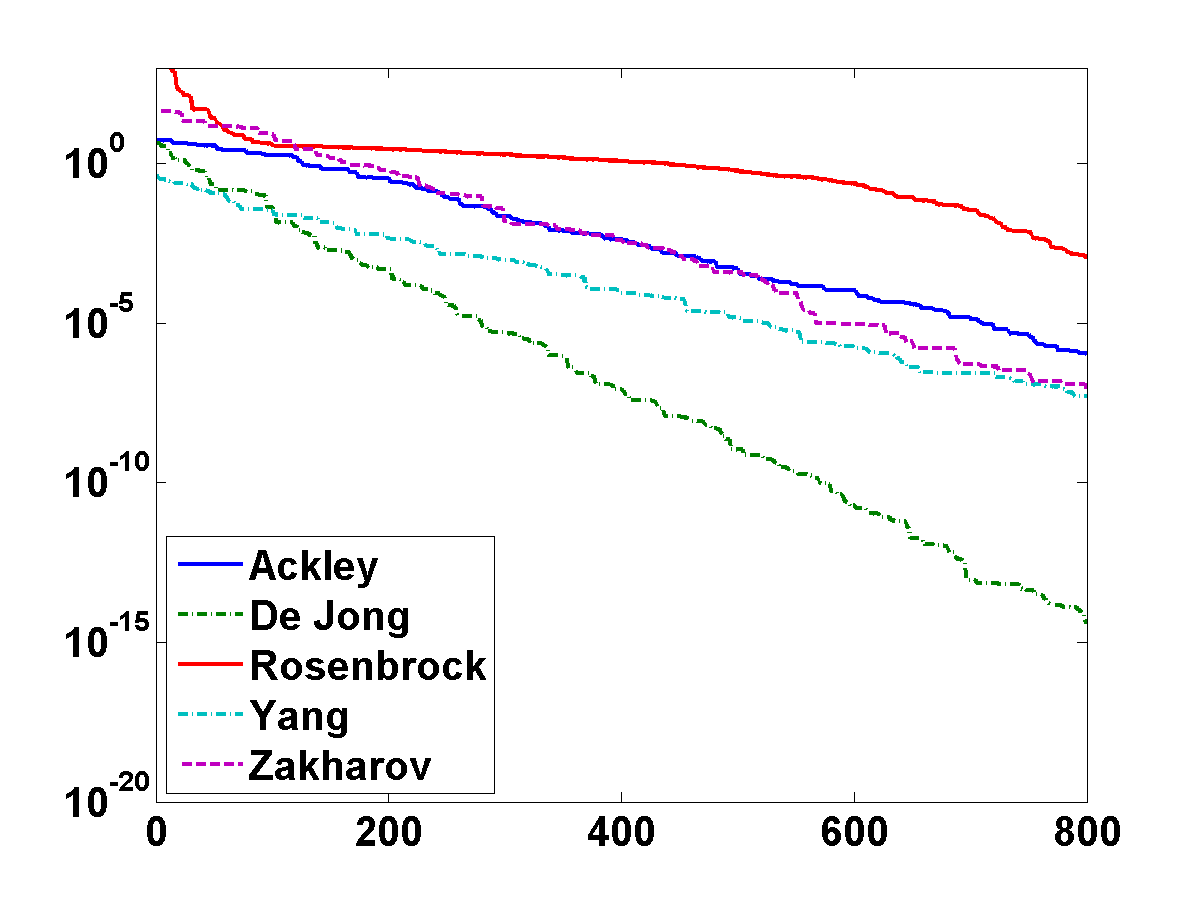}}
\caption{Convergence of five test functions using the flower algorithm.}
\label{fig-fiveplot}
\end{figure}

By using the flower algorithm with $n=20$, $\lam=1.5$, $p=0.8$
and a fixed number of iterations $t=1000$, we can find the global minima
for all the above 5 functions for $d=4$. The convergence graphs for all these functions
are summarized and shown in Fig.~\ref{fig-fiveplot}. As we can see,
they all converge quickly in an almost exponential manner, except for
Rosenbrock's function which has a narrow valley. Once the search has gone through
some part of the valley during iterations, its convergence becomes exponentially
with a higher slope, though the rate of convergence is still lower compared with those for other functions.

Though the theoretical analysis proves that FPA will converge, it is worth pointing out the
the rate of convergence is still influenced by both the algorithmic structure and its parameter settings. The convergence analysis does not provide much information about how quickly the algorithm may converge for a given problem, and consequently parameter tuning may be needed in practice to find the best parameter settings to give a higher convergence rate.

\section{Conclusions}

The flower pollination algorithm  is an efficient optimization algorithm with a wide range of applications. We have provided the first results on the convergence analysis of this algorithm.
By using the Markov models, we have proved that the flower pollination algorithm has guaranteed global convergence,
which laid the theoretical foundation for this algorithm and showed why it is efficient
in applications.  Then, we have used a set of five different functions with diverse properties to
show that FPA can indeed converge very quickly in practice.

It is worth pointing out that the current results are mainly for the standard flower pollination algorithm. It will be useful if further research can focus on the extension of the proposed methodology to analyze
the convergence of the full flower pollination algorithm and its variants.
Ultimately, it can be expected that
the proposed method can be used to analyze other metaheuristic algorithms as well.

\section*{References}
\bibliographystyle{elsarticle-num}
\bibliography{<your-bib-database>}

%% End of text %%

\end{document}